\documentclass[journal]{IEEEtran}

\usepackage{cite}

\usepackage{graphicx}

\usepackage[cmex10]{amsmath}

\usepackage{array}

\usepackage{url}
  \usepackage{subfigure}
\usepackage{makecell,multirow,diagbox}
\newcommand{\brem}{\medskip \begin{remark}  }
\newcommand{\erem}{\hfill \rule{1mm}{2mm}\medskip
\end{remark} }
\newtheorem{remark}{\bf Remark}[section]
\begin{document}

\title{Granular Optimal Load-Side Control of Power Systems with Electric Spring Aggregators}

\author{Congchong~Zhang,~\IEEEmembership{Student Member,~IEEE,}
        Tao~Liu,~\IEEEmembership{Member,~IEEE,}
        and~David~J.~Hill,~\IEEEmembership{Life~Fellow,~IEEE}
        \vspace{-14.4pt}
\thanks{The work described in this paper was fully supported by a grant from the Research Grants Council of the Hong Kong Special Administrative Region under Theme-based Research Scheme through Project No. T23-701/14-N.}
}

\markboth{Journal of \LaTeX\ Class Files,~Vol.~13, No.~9, September~2014}%
{Congchong \MakeLowercase{\textit{et al.}}: Bare Demo of IEEEtran.cls for Journals}



\maketitle

\begin{abstract}

To implement controllable loads for frequency regulation in transmission networks in a practical way, the control scheme needs to be granulated down at least to subtransmission networks since loads in transmission networks are usually the aggregation of lower voltage networks. However, not only frequency but also bus voltage will be affected by active power changes in subtransmission networks due to a higher R/X ratio of transmission lines. Further, the costs for loads participating in frequency and voltage regulation should also be considered. In this paper, a control scheme is proposed for electric spring (ES) aggregators which consist of back-to-back ESs and exponential type of noncritical loads in subtransmission networks. A distributed optimization which aims to minimize the costs and implements both frequency and voltage regulation is adopted for ES aggregators to obtain new active and reactive power setpoints by sharing information with neighbors. Power consumption of each ES aggregator is then adjusted accordingly to conduct frequency and voltage regulation simultaneously. Simulation results show that ES aggregators are able to achieve required active power response and regulate frequency cooperatively, and meanwhile maintain bus voltages within the acceptable range with minimized costs under the proposed control scheme.

\end{abstract}


\IEEEpeerreviewmaketitle

\section{Introduction}

Traditional frequency control is implemented on the generation side by adjusting the mechanical power inputs of generators to follow demand, which consists of three control layers operating at different timescales \cite{kundur1994power}, i.e., droop control (primary frequency control), automatic generation control (secondary frequency control) and economic dispatch (tertiary frequency control). However, due to the uncertainty and intermittency of renewable power, this traditional control method may be inadequate to keep the system frequency at its nominal value in the future \cite{windenergy_2005}. The situation may become even worse as the penetration level of renewable power increases. To integrate more renewables, a large extra quantity of spinning reserves will be required \cite{xu2011demand}, which may lead to much higher operation costs. 

To overcome these issues, load-side frequency control (or demand response) has drawn a lot of attention lately because of its advantages such as instantaneous response, potentially lower costs and highly distributed availability throughout the grid \cite{callaway2011achieving,xu2011demand}. Different control methods have been proposed to automatically adjust power consumption of loads in transmission networks for frequency regulation \cite{tao_nondisrupt,ZHANG_GM,xu2011demand, zhao2015distributed,zhangpscc, mallada2014distributed,Zhengyu_GM,callaway2011achieving}. In our previous work \cite{tao_nondisrupt,ZHANG_GM}, a switched consensus-based distributed control method has been proposed for controllable loads in transmission networks which can work in two different modes, i.e., the frequency regulation mode (FRM) and load recovery mode (LRM), to achieve both being fully responsive and non-disruptive when participating in frequency control. 

However, due to the hierarchical structure of electric power grids with respect to different voltage levels \cite{hierarchy}, all aforementioned load-side control methods for transmission networks need to be granulated down to subtransmission networks to be implemented in practical ways since loads in transmission networks are aggregates of large numbers of physical or aggregate loads in subtransmission networks. Further, in subtransmission networks the active power changes will affect not only the frequency but also bus voltages due to a higher $R/X$ ratio of lines and cables. Moreover, the costs for loads providing active and reactive power supports need to be considered. Therefore, an optimization algorithm which can minimize the costs and meanwhile takes both frequency and voltage regulation into consideration is required.

In our previous work \cite{zhangccpowertech}, the electric spring (ES) aggregators that consist of large numbers of ESs and noncritical loads are adopted to achieve the required active power response and regulate bus voltages simultaneously. However, it does not consider the costs of ES aggregators for providing a required power response. Further, the ESs used in \cite{zhangccpowertech} are an old version of ES (ES-2) which can only support limited active power and reactive power within a certain range due to the capacity limit of the battery and characteristics of the noncritical load. Moreover, the noncritical loads that cascade with ESs are usually considered as constant impedance or constant resistive loads in most of the existing works \cite{zhangccpowertech,b2b_voltageES,yan_esb2b,hui2012electric}. However, many other load characteristics are used in practice such as constant current, constant power and exponential types in power grids. 

Thus, in this paper a distributed optimization algorithm is proposed for ES aggregators in subtransmission networks to provide the required active power response and meanwhile maintain bus voltages with minimized costs. Due to the nonlinearity and low computational efficiency of the traditional AC power flow calculation, a decoupled linear power flow model proposed in \cite{DLPF} is used to minimize the costs of ES aggregators. In the optimization algorithm, the required active power response of each aggregate controllable load (i.e., the control output of each load-side controller in transmission networks proposed in \cite{ZHANG_GM}) and bus voltages in subtransmission networks are considered as local constraints for each ES aggregator. By sharing information with neighboring aggregators, each ES aggregator can accomplish the corresponding power flow calculation and obtain an optimal solution (the active and reactive powers setpoints) in a cooperative way. Then, the ES voltage will be adjusted according to the obtained setpoints and the power consumption of noncritical loads will be changed subsequently to regulate the frequency and voltage simultaneously with minimized costs. Differently from using ES-2 in \cite{zhangccpowertech}, in this paper the third version of the ES is adopted which consists of two half-bridge inverters configured in a back-to-back structure (ES-B2B). This new type of the ES has a larger power support capability by replacing the battery with an ac-to-dc shunt inverter \cite{b2b_voltageES,yan_esb2b}. Moreover, an exponential load model which can represent a combination of different types of loads is adopted for the noncritical load in this paper.

The rest of this paper is organized as follows. In Section \uppercase\expandafter{\romannumeral2}, the transmission network model with aggregate loads and both generation-side and load-side frequency control are introduced. Section \uppercase\expandafter{\romannumeral3} introduces the subtransmission network and ES aggregator model. Section \uppercase\expandafter{\romannumeral4} illustrates the control scheme for ES aggregators. Case studies are presented and analyzed in Section \uppercase\expandafter{\romannumeral5}. The paper ends by conclusions in Section \uppercase\expandafter{\romannumeral6}.

\section{Transmission Network Model and Frequency Control Review}\label{SC}

In this paper, the following standard assumptions are made for the transmission network:
 
1) The transmission network is connected in which transmission lines are lossless and characterized by reactances $x_{ij}=x_{ji}$. 

2) The magnitude of voltage $|V_{i}|$ of each bus in the transmission network is fixed, then the active power flows between buses will not be affected by the voltages. 


The structure-preserving model proposed in \cite{bergen1981structure} with aggregate loads is adopted to model the transmission network,

\begin{subequations}\label{eq1}
\begin{align}
 \dot{\delta_{i}}&=\omega_{i},~~~~~~~~~~~~~~~~~~~~~~~~~~~~~~~~~~~~~~~~i\in\mathcal{G}\\
M_{i}\dot{\omega}_{i}&=-D_{i}\omega_{i}+P_{m_{i}}-\sum_{j=1}^{N}b_{ij}\sin(\delta_{i}-\delta_{j}),i\in\mathcal{G}\label{eq1b}\\
 D_{i}\dot{\delta_{i}}&=u_{i}-\sum_{j=1}^{N}b_{ij}\sin(\delta_{i}-\delta_{j})-P_{D_{i}},~~~~~~i\in\mathcal{L}\label{eq1c}
\end{align}
\end{subequations}
where $\mathcal{N}= \{1, 2, ..., N\}= \mathcal{G}\bigcup\mathcal{L}$ is the index set of all the buses. The index sets of generator buses and load buses are denoted by $\mathcal{G}=\{1, 2, ..., N_{G}\}$ and $\mathcal{L}=\{1, 2, ..., N_{L}\}$ with cardinalities $N_{G}$ and $N_{L}$, respectively.

For all $i\in\mathcal{N}$, the coefficient $b_{ij}=\frac{|V_{i}||V_{j}|}{x_{ij}}$ is acquired based on assumption 1) and 2). For each generator $i\in\mathcal{G}$, parameters and variables $M_{i}, D_{i},\delta_{i}, \omega_{i}$ and $P_{m_{i}}$ represent the inertia constant, damping coefficient, power angle, power angular velocity and mechanical power input, respectively. For each load $i\in\mathcal{L}$, symbols $D_{i},\delta_{i}, u_{i}$ and $P_{D_{i}}$ represent the frequency-dependence coefficient, bus voltage phase angle, actual response of the aggregate controllable load with capacity limits $\underline{u}_{i}\le u_{i} \le \overline{u}_{i}$ and uncontrollable constant load, respectively.

In our previous work \cite{tao_nondisrupt,ZHANG_GM}, we proposed a cooperative control method to coordinate generation-side control and load-side control. On the generation side, some of the generators are selected to participate in AGC and adjust their setpoints every few seconds to restore the frequency to the nominal value, whereas other generators have droop control only \cite{kundur1994power}. On the load side, a switched consensus-based distributed control method has been proposed for the load-side controller. It works in the FRM to restore the system frequency after disturbances, and then switches to the LRM to recover aggregate controllable loads to their nominal values after the frequency goes back to an acceptable region. Thus, for each load $i\in\mathcal{L}$, the distributed controller is proposed as follows \cite{ZHANG_GM}, 
\begin{subequations}\label{eq2}
\begin{align}
\text{FRM:}~ \dot{r}_{i}(t)&=\sum_{j=1}^{N}a_{ij}(r_{j}(t)-r_{i}(t)), t \in [t_{m},t_{m}^{'}+\tau)\\
 r_{i}(t^{+})&=d_{i}(t), t=t_{m}, m=1,2,...\label{ri}\\
 \bar{u}_{i}(t)&=\!-\!K_{L_{i}}^{I}r_{i}(t)\!-\!K_{L_{i}}^{P}\Delta\omega_{i}(t), t\!\in\![t_{m},t_{m}^{'}+\tau)\label{ui}\\
\text{LRM:}~ \dot{r}_{i}(t)&=-\alpha_{i} r_{i}(t), t\in [t_{m}^{'}+\tau,t_{m+1})\\
\bar{u}_{i}(t)&=-K_{L_{i}}^{I}r_{i}(t), t\in[t_{m}^{'}+\tau,t_{m+1})\label{loadrecovery}
 \end{align}
 \end{subequations}
where parameters $K_{L_{i}}^{I}>0$ and $K_{L_{i}}^{P}>0$ represent the integral and proportional control gains, respectively. The variable $\Delta\omega_{i}=\omega_{i}-\omega_{s}$ is the angular velocity deviation where $\omega_{s}=2\pi f_{s}$ is the nominal angular velocity of the system with $f_{s}$ representing the nominal value of the system frequency. The variable $\bar{u}_{i}$ is the required response of each aggregate controllable load. The time instant at which the $m$th time the system frequency exceeds the pre-set frequency band due to an active power mismatch is denoted by $t_{m}$. The power imbalance of each bus at $t=t_{m}$ denoted by $d_{i}$ is defined as: $d_i=P_{m_{i}}-\sum_{j=1}^{N}b_{ij}\sin(\delta_{i}-\delta_{j})$ for each bus $i\in\mathcal{G}$ and $d_i=u_{i}-\sum_{j=1}^{N}b_{ij}\sin(\delta_{i}-\delta_{j})-P_{D_{i}}$ for each bus $i\in\mathcal{L}$. The parameter $\alpha_{i}$ is the load recovery rate. The matrix $A=(a_{ij})\in\mathcal{R}^{N\times N}$ represents the topology of the communication network for the controllers in transmission networks. The time instant when the frequency recovers to the acceptable region is denoted by $t_{m}^{'}$, and the constant $\tau$ is the dwell-time aiming to reduce unnecessary switching between the FRM and LRM. It should be noted that the generator buses are also included in the distributed consensus algorithm, for helping load-side controllers to acquire the average power imbalance of the system \cite{ZHANG_GM}.

Let $\phi(t)$ denote the switching signal. Further, let $\phi(t)=1$ represent the aggregate controllable loads working in the LRM, and $\phi(t)=0$ represent the aggregate controllable loads  working in the FRM. At each $t=t_{m}$, the control centre sends $\phi(t)=0$ to inform each controller. Then each controller reinitializes its state $r_{i}(t)$ according to (\ref{ri}) and starts to discover the average power imbalance of the system by communicating with its neighbors. The control output $\bar{u}_{i}$ of each load-side controller is then adjusted accordingly to help generators to restore frequency. When the frequency goes back to the acceptable region for $\tau$ s, each load-side controller will begin to recover aggregate controllable loads to their nominal values when it receives the signal $\phi(t)=1$ from the control centre at $t=t_{m}^{'}+\tau$ s, and the generators will gradually take full responsibility for the active power imbalance in the system. Details of the proposed control method can be found in \cite{ZHANG_GM}.

It should be noted that the variable $\bar{u}_{i}$ in (\ref{eq2}) is the required but not the acutal response of each aggregate controllable load which needs to be sent to subtransmission networks where it can be implemented by adjusting physical loads coordinately. Thus, a coordinative control method should be proposed by which the actual aggregate response of controllable loads in subtransmission networks i.e., $u_{i}$ in (\ref{eq1c}) can track the corresponding required response $\bar{u}_{i}$, and will be discussed in Section \ref{Implement}.

\section{Subtransmission Network Model and ES Aggregator Modeling}\label{subtrans}

Loads in the transmission network model are usually aggregated by large numbers of physical loads in subtransmission networks, and hence the load-side control signals obtained in the transmission network need to be sent to subtransmission networks to be implemented in a practical way. Thus, each load bus in the transmission network is assumed to be the aggregation of a subtransmission network in this paper. Differently from transmission networks, transmission lines in subtransmission networks have a higher $R/X$ ratio. Consequently, not only the frequency but also bus voltages will be affected by active power changes, and hence both of them need to be regulated simultaneously when any active power mismatches occur. 

In this situation, a smart load consisting of an electric spring (ES-B2B) and its associated noncritical load can be adopted to deal with this issue effectively since it can provide active power and reactive power support simultaneously \cite{yan_esb2b}. Therefore, we use an ES aggregator and a critical load connected in parallel at each load bus to represent the load in subtransmission networks. Two cascaded components are included in each ES aggregator: a large control capacity ES-B2B and an aggregate noncritical load. Fig. \ref{ES figure} shows the configuration of such a load in which the symbols $\mathbf{V}_{s},\mathbf{V}_{es},\mathbf{V}_{nl}$ and $\mathbf{I}$ represent the bus voltage, ES-B2B voltage, noncritical load voltage and current, respectively.
The structure of the ES-B2B and the load model will be discussed in detail in the following subsections.

\subsection{Basic Concept of ES}\label{basic concept}

\begin{figure}[h]
\centering
\includegraphics[width=2.7in]{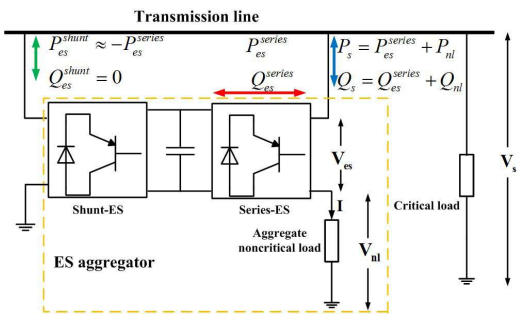}
\caption{Configuration of a load with an ES aggregator and a critical load}
\label{ES figure}
\end{figure}

So far, three versions of ESs have been conceived. The third version of the ES (i.e., ES-B2B) consisting of two half-bridge inverters is adopted in this paper \cite{yan_esb2b}. The configuration of the ES-B2B is shown in Fig. \ref{ES figure}. Similar to the first two versions, one inverter of the ES-B2B (Series-ES) is connected in series to a noncritical load. Unlike using a capacitor in ES-1 or a battery in ES-2 as energy storage, the ES-B2B uses a bidirectional ac-to-dc inverter (Shunt-ES) with the ac side connected to the power grid (see Fig. \ref{ES figure} where the Shunt-ES works in parallel with the Series-ES plus aggregate noncritical load). The active power supplied (or consumed) by the Series-ES is equal to the active power consumed (or supplied) by the Shunt-ES from the grid, and the Shunt-ES does not consume any reactive power, i.e., $P_{es}^{shunt}=-P_{es}^{series}$ and $Q_{es}^{shunt}=0$ \cite{yan_esb2b}. 
Thus, the total active power $P_{sl}$ and reactive power $Q_{sl}$ of an ES aggregator can be represented as follows,
\begin{subequations}\label{PQtotal}
\begin{align}
P_{sl}&=P_{s}+P_{es}^{shunt}=P_{nl}\label{Ptotal}\\
Q_{sl}&=Q_{s}+Q_{es}^{shunt}=Q_{nl}+Q_{es}^{series}\label{Qtotal}
 \end{align}
 \end{subequations}
where the symbols $P_{s}$ and $Q_{s}$ represent the total active power and reactive power of the series-ES and associated noncritical load, respectively; and $P_{nl}$ and $Q_{nl}$ are the active and reactive power consumption of the noncritical load, respectively.

\subsection{Load Model}
The operation of the ES-B2B affects the voltage of the noncritical load significantly and consequently its power consumption as well. So far, only constant resistive or impedance loads have mostly been used as the noncritical load in the existing ES studies (e.g., \cite{yan_esb2b,hui2012electric,b2b_voltageES,zhangccpowertech}). To extend the adoption of ESs, a general load model should be considered to verify the proposed method with ESs \cite{zheng2017critical}. In this paper, the noncritical load is modeled as an exponential load whose active power and reactive power can be expressed by the following exponential functions of voltage, respectively
\begin{subequations}\label{expo}
\begin{align}
P_{nl}=P_{0}\left(\frac{V_{nl}}{V_{0}}\right)^{\alpha_{p}}\label{expoP}\\
Q_{nl}=Q_{0}\left(\frac{V_{nl}}{V_{0}}\right)^{\alpha_{q}}\label{expoQ}
 \end{align}
 \end{subequations}
where the constants $P_{0},~Q_{0}$ and $V_{0}$ are the nominal active power, reactive power and magnitude of the noncritical load voltage, respectively. For simplicity, $V_{0}$ of all load buses are assumed to be one unit in this paper. The variable $V_{nl}$ is the voltage magnitude on the noncritical load with limits $0<\underline{V}_{nl}\le V_{nl}\le \overline{V}_{nl}$, and $\alpha_{p}$ ($\alpha_{p}\neq0$) and $\alpha_{q}$ are the exponential coefficients of the noncritical load. By substituting (\ref{expoP}) into (\ref{expoQ}), we can rewrite (\ref{expoQ}) as follows,
\begin{equation}\label{Qcl}
Q_{nl}=Q_{0}\left(\frac{P_{nl}}{P_{0}}\right)^{\frac{\alpha_{q}}{\alpha_{p}}}.
\end{equation}
The method on how to control ES-B2Bs to regulate frequency and voltage will be discussed in next section. 

\brem
We assume that the critical loads are constant power in this paper. Because the bus voltages will be maintained within small ranges around the nominal values due to the control actions of ES aggregators, the power consumption of critical loads will consequently be maintained around their nominal value. Thus, for simplicity we ignore these small power changes and assume the critical loads are unchanged.
\erem

\brem
Although the ES aggregator has a relatively large control capacity as it consists of many smart loads, a single ES aggregator may still not fulfill the requirement of load-side control in transmission networks. Then, the problem how to coordinate a number of ES aggregators in a subtransmission network such that their aggregate control actions can satisfy the corresponding requirement arises. This question will be answered in the next section.
\erem

\brem
Regarding the frequency-dependence load in (\ref{eq1c}), we assume that the frequency-dependence coefficient $D_i$ of each load $i$ is small, and for simplicity, we still use the term $D_i\dot{\delta_i}$ in the transmission network. However, the method on how to granulate this term down to the subtransmission network is of importance and deserves attention.
\erem

\section{Optimal Control for ES Aggregator}\label{Implement}

In subtransmission networks, both the frequency and bus voltages will be affected by active power changes due to the higher $R/X$ ratio of transmission lines, and hence both need to be taken into account by ES aggregators when active power mismatch occurs. However, bus voltages do not have to all hold to the nominal value and are allowed to vary within required limits \cite{kundur1994power}. Further, the amount of active power and reactive power that each ES aggregator should provide needs to be coordinated and optimized to minimize the costs and fulfill the active power response (transmission level) and bus voltage limits requirements in the meantime. To solve these problems, we propose a control method and a distributed optimization algorithm. The designed control method is used to adjust active and reactive powers of ES aggregators according to the reference setpoints obtained from the proposed optimization algorithm, which will be discussed in this section.

\subsection{Active and Reactive Power Control of ES Aggregators}
To achieve the required active power, i.e., $\bar{u}_{i}$ at the point of common coupling (PCC), and meanwhile keep bus voltages within required limits, both the active power and reactive power of each ES aggregator need to be adjusted simultaneously. Therefore, a $d-q$ transformation is adopted to decouple the ES-B2B voltage $\mathbf{V}_{es}$, noncritical load voltage $\mathbf{V}_{nl}$ and current $\mathbf{I}$ into $d$ and $q$ components, respectively, as follows
\begin{subequations}\label{dq}
\begin{align}
\mathbf{V}_{es}&=V_{es\_{d}}+jV_{es\_{q}}\\
\mathbf{V}_{nl}&=V_{nl\_{d}}+jV_{nl\_{q}}\label{Vesdq}\\
\mathbf{I}&=I_{d}+jI_{q}.
\end{align}
\end{subequations}
More details of the $d-q$ transformation can be found in \cite{zhangccpowertech}.

%

Let $\mathcal{N}_{i}=\{1,2,...,N_{i}\}$ and $\mathcal{N}_{es}^{i}=\{1,2,...,N_{es}^{i}\}$ denote the index set of all the buses and buses with ES aggregators in each subtransmission network $i\in\mathcal{L}$, respectively. 
For each ES aggregator $j\in\mathcal{N}_{es}^{i}$, the bus voltage is used as the reference, i.e., $\mathbf{V}_{s_{j}}=V_{s_{j}}\angle0^{\circ}$, to obtain the corresponding ES voltage. Easily, we can acquire
\begin{subequations}\label{Vnl}
\begin{align}
V_{nl\_d_{j}}+V_{es\_d_{j}}&=V_{s_{j}}\label{Vnld}\\
V_{nl\_q_{j}}+V_{es\_q_{j}}&=0\label{Vnlq}.
\end{align}
\end{subequations}
From (\ref{expoP}) and (\ref{Vesdq}), we have
 \begin{equation}\label{sqrt}
 \sqrt{V_{nl\_d_{j}}^{2}+V_{nl\_q_{j}}^{2}}=V_{0_{j}}\left(\frac{P_{nl_{j}}}{P_{0_{j}}}\right)^{\frac{1}{\alpha_{p_{j}}}}.
 \end{equation}
By using (\ref{PQtotal}) and (\ref{Qcl}), the active and reactive power consumption $P_{nl_{j}}$ and $Q_{nl_{j}}$ of the noncritical load and reactive power consumption $Q_{sl_{j}}$ of the ES aggregator can be represented as
 \begin{subequations}\label{PQnl}
\begin{align}
P_{nl_{j}}&=V_{nl\_d_{j}}I_{d_{j}}+V_{nl\_q_{j}}I_{q_{j}}=P_{j}^{*}\label{Pnl}\\
Q_{nl_{j}}&=V_{nl\_q_{j}}I_{d_{j}}-V_{nl\_d_{j}}I_{q_{j}}=Q_{0_{j}}\left(\frac{P_{nl_{j}}}{P_{0_{j}}}\right)^{\frac{\alpha_{q_{j}}}{\alpha_{p_{j}}}}\label{Qnl}\\
 Q_{sl_{j}}&=-V_{s_{j}}I_{q_{j}}=Q_{j}^{*}\label{Qsl}
\end{align}
\end{subequations}
where the required active power and reactive power of each ES aggregator $j$ are denoted by $P_{j}^{*}$ and $Q_{j}^{*}$, respectively. Here we drop the $sl$ subscript for simplicity.

By using (\ref{Vnl}), (\ref{sqrt}) and (\ref{PQnl}), the $d$ and $q$ components of each ES-B2B voltage setpoints with saturation limits $\underline{V}_{es\_d}\le V^{*}_{es\_d_{j}} \le \overline{V}_{es\_d}$ and $\underline{V}_{es\_q}\le V^{*}_{es\_q_{j}} \le \overline{V}_{es\_q}$ are acquired as follows 
 \begin{equation}\label{apcm_2}
\begin{split}
V^{*}_{es\_d_{j}}&\!=\!V^{*}_{s_{j}}\!-\!\frac{Q_{0_{j}}\left(\frac{P^{*}_{j}}{P_{0_{j}}}\right)^{\frac{\alpha_{q_{j}}+2}{\alpha_{p_{j}}}}Q^{*}_{j}+\left(\frac{P^{*}_{j}}{P_{0_{j}}}\right)^{\frac{1}{\alpha_{p_{j}}}}P^{*}_{j}\sqrt{\Delta}}{P^{*2}_{j}+Q_{0_{j}}^{2}\left(\frac{P^{*}_{j}}{P_{0_{j}}}\right)^{\frac{2\alpha_{q_{j}}}{\alpha_{p_{j}}}}}V_{0_{j}}\\
V^{*}_{es\_q_{j}}&=\frac{\left(\frac{P^{*}_{j}}{P_{0_{j}}}\right)^{\frac{2}{\alpha_{p_{j}}}}P^{*}_{j}Q^{*}_{j}-Q_{0_{j}}\left(\frac{P^{*}_{j}}{P_{0_{j}}}\right)^{\frac{\alpha_{q_{j}}+1}{\alpha_{p_{j}}}}\sqrt{\Delta}}{P^{*2}_{j}+Q_{0_{j}}^{2}\left(\frac{P^{*}_{j}}{P_{0_{j}}}\right)^{\frac{2\alpha_{q_{j}}}{\alpha_{p_{j}}}}}V_{0_{j}}
\end{split}
\end{equation}
where $V^{*}_{s_{j}}$ is the required bus voltage and $\Delta=P^{*2}_{j}-\left(\frac{P^{*}_{j}}{P_{0_{j}}}\right)^{\frac{2}{\alpha_{p_{j}}}}Q^{*2}_{j}+Q_{0_{j}}^{2}\left(\frac{P^{*}_{j}}{P_{0_{j}}}\right)^{\frac{2\alpha_{q_{j}}}{\alpha_{p_{j}}}}$.  Two closed-loop PI controllers are adopted to enable the actual $d$ and $q$ components of each ES-B2B voltage $V_{es\_d_{j}}$ and $V_{es\_q_{j}}$ to hold to the setpoints $V^{*}_{es\_d_{j}}$ and $V^{*}_{es\_q_{j}}$, respectively. The setpoints $V^{*}_{es\_d_{j}}$ and $V^{*}_{es\_q_{j}}$ are determined by the required power consumption of the ES aggregator, i.e., $P^{*}_{j}$ and $Q^{*}_{j}$ obtained from the proposed distributed optimization algorithm which will be discussed in the next subsection.

\subsection{Distributed Optimization Over ES Aggregators}
As discussed in the last subsection, to regulate the frequency and voltage simultaneously, new active and reactive power setpoints are required by each ES aggregator to adjust its ES-B2B voltage as given in (\ref{apcm_2}). Moreover, the costs for the active and reactive power support by ES aggregators should also be considered. Therefore, optimal power flow calculations aiming to minimize the costs and from which each ES aggregator can obtain its corresponding active and reactive power setpoints need to be conducted. Since there may be large numbers of ES aggregators scattered in subtransmission networks, centralized optimization may be infeasible to deal with this problem. In contrast, distributed optimization can address this problem effectively in which each ES aggregator only needs to communicate with its neighbors to acquire the optimal solution cooperatively. However, due to the nonlinearity and low computational efficiency, the traditional AC power flow model may be infeasible to be adopted in the distributed optimization, and hence a linear power flow model is used.

In this paper, a decoupled linearized power flow (DLPF) model proposed in \cite{DLPF} is adopted to approximate the AC power flow model. We illustrate the DLPF model for one subtransmission network as an example, which will then be adopted to all subtransmission networks in this paper, and hence the subscript $i$ standing for the $i$th subtransmission network is dropped for simplicity. Thus, the matrix form of this DLPF model is given as follows,

\begin{equation}\label{DLPF1}
\begin{split}
\begin{bmatrix}
\textbf{P} \\
\textbf{Q}
\end{bmatrix}
=-
\begin{bmatrix}
-\textbf{G} & \textbf{B}^{'} \\
\textbf{B} & \textbf{G}
\end{bmatrix}
\begin{bmatrix}
 \textbf{V}\\
\boldsymbol{\theta}
\end{bmatrix}
\end{split}
\end{equation}
%
where the vectors $\textbf{P}=(P_{1},P_{2},...,P_{N_{i}})^{T},~\textbf{Q}=(Q_{1},Q_{2},...,Q_{N_{i}})^{T},~ \textbf{V}=(V_{1},V_{2},...,V_{N_{i}})^{T}$ and $\boldsymbol{\theta}=(\theta_{1},~\theta_{2},...,\theta_{N_{i}})^{T}$ are the bus injected active powers, reactive powers, bus voltage magnitudes and phase angles, respectively. The matrices $\textbf{G}\in\mathcal{R}^{N_{i}\times N_{i}}$ and $\textbf{B}\in\mathcal{R}^{N_{i}\times N_{i}}$ are the real part and imaginary part of the admittance matrix of the $i$th subtransmission network, respectively, and $\textbf{B}^{'}\in\mathcal{R}^{N_{i}\times N_{i}}$ is the imaginary part of admittance matrix without shunt elements. Define the vector $\textbf{x}=(\textbf{P}^{T},\textbf{Q}^{T},\textbf{V}^{T},\boldsymbol{\theta}^{T})^{T}\in\mathcal{R}^{4N_{i}}$. Hence, the model (\ref{DLPF1}) can be rewritten as
\begin{equation}\label{DLPF2}
\begin{split}
\underbrace{\begin{bmatrix}
\textbf{I} & \textbf{O} & -\textbf{G} & \textbf{B}^{'} \\
\textbf{O} & \textbf{I} & \textbf{B} & \textbf{G}
\end{bmatrix}}_{\textbf{W}}
\textbf{x}
=
\boldsymbol{0}.
\end{split}
\end{equation}
where the matrices $\textbf{I}\in\mathcal{R}^{N_{i}\times N_{i}}$ and $\textbf{O}\in\mathcal{R}^{N_{i}\times N_{i}}$ are the identity matrix and zero matrix, respectively, and $\boldsymbol{0}\in\mathcal{R}^{2N_{i}}$ is a zero vector. More details of the DLPF model can be found in \cite{DLPF}.

Thus, for each subtransmission network, the optimization problem can be described as follows,
\begin{equation}\label{centraloptimal}
\text{minimize}~~f^{total}=\sum^{N_{es}^{i}}_{j=1}f_{j}(\textbf{x})
\end{equation}
subject to
\begin{equation}
\label{centraloptimal_subject}
\textbf{W}~\textbf{x}=\boldsymbol{0},~\textbf{x}\in\Omega
\end{equation}
where the cost function of the power compensation by each aggregator $j$ is assumed to be convex and selected to be a quadratic function, namely $f_{j}(\textbf{x})\!=\!h_{j}(P_{j}\!-\!P_{0_{j}})^{2}\!+\!g_{j}(Q_{j}\!-\!Q_{0_{j}})^{2}$ where $h_{j}$ and $g_{j}$ are the cost coefficients of active and reactive powers, respectively. We assume that the cost function $f_{j}$ is only known by aggregator $j$ for privacy issues. The set $\Omega_{i}=[\underline{\textbf{X}},\overline{\textbf{X}}]$ is the global constraint with $\underline{\textbf{X}}\!=\!(\underline{\textbf{P}}^{T},\underline{\textbf{Q}}^{T},\underline{\textbf{V}}^{T},\underline{\boldsymbol{\theta}}^{T})^{T}$ and $\overline{\textbf{X}}\!=\!(\overline{\textbf{P}}^{T},\overline{\textbf{Q}}^{T},\overline{\textbf{V}}^{T},\overline{\boldsymbol{\theta}}^{T})^{T}$, respectively, in which $\underline{\textbf{P}}\!=\!( \underline{P}_{1},\underline{P}_{2},...,\underline{P}_{N_{i}})^{T}$, $\overline{\textbf{P}}\!=\!( \overline{P}_{1},\overline{P}_{2},...\overline{P}_{N_{i}})^{T}$, $\underline{\textbf{Q}}\!=\!( \underline{Q}_{1},\underline{Q}_{2},...,\underline{Q}_{N_{i}})^{T}$, $\overline{\textbf{Q}}\!=\!( \overline{Q}_{1},\overline{Q}_{2},...\overline{Q}_{N_{i}})^{T}$, $\underline{\textbf{V}}\!=\!(\underline{V}_{s_{1}},\underline{V}_{s_{2}},...,\underline{V}_{s_{N_{i}}})^{T}$, $\overline{\textbf{V}}\!=\!(\overline{V}_{s_{1}},\overline{V}_{s_{2}},...,\overline{V}_{s_{N_{i}}})^{T}$, $\underline{\boldsymbol{\theta}}\!=\!(\underline{\theta}_{1},\underline{\theta}_{2},...,\underline{\theta}_{N_{i}})^{T}$ and $\overline{\boldsymbol{\theta}}\!=\!(\overline{\theta}_{1},\overline{\theta}_{2},...,\overline{\theta}_{N_{i}})^{T}$ are the lower and upper limits of active powers, reactive powers, bus voltage magnitudes and phase angles, respectively. 

\begin{figure}[htbp]
\centering
\includegraphics[width=2.3in]{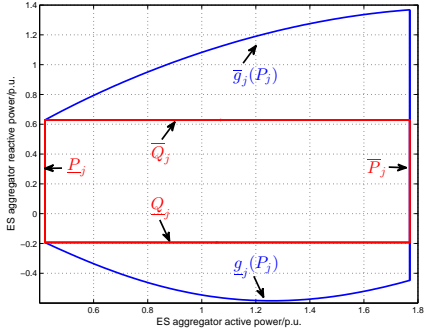}
\caption{The relationship between active and reactive powers of ES aggregator}
\label{PQrelationship}
\end{figure}

It should be noted that bus $1$ in each subtransmission network is assumed to be the PCC, and we set $\underline{V}_{s_{1}}=\overline{V}_{s_{1}}=V_{P}$, $\underline{\theta}_{1}=\overline{\theta}_{1}=\theta_{P}$, $\underline{Q}_{1}=\underline{Q}_{P}$ and $\overline{Q}_{1}=\overline{Q}_{P}$.  We assume that a PMU is installed at the PCC in each subtransmission network such that the active power flow $P_{pcc}$ at the PCC can be measured. Transmission losses $P_{loss}$ of each subtransmission network need to be considered because of the high $R/X$ ratio of lines in substransmission networks. Therefore, we further assume there is linear relationship between the transmission losses $P_{loss}$ and active power flow $P_{pcc}$, i.e., $P_{loss}=D^{l}P_{pcc}+C$ where $D^{l}$ is the coefficient and $C$ is a constant \cite{wu2016line}.
Thus, to achieve the required active power, the active power limits of bus $1$ in each subtransmission network is set to be equal to the control reference signal $\bar{u}$ from the load-side controller at the corresponding bus in the transmission network plus transmission losses of the subtransmission networks, i.e., $\overline{P}_{1}=\underline{P}_{1}=\bar{u}+P_{loss}$.  For buses $j\in\mathcal{N}_{i}\setminus\{1\}$, the limits of magnitudes and phase angles of bus voltages are uniformly set as $\underline{V}_{s_{j}}=\underline{V}_{s}$, $\overline{V}_{s_{j}}=\overline{V}_{s}$, $\underline{\theta}_{j}=\underline{\theta}$ and $\overline{\theta}_{j}=\overline{\theta}$, to guarantee bus voltages not to exceed the required limits. For the active power limits, from (\ref{expoP}) we have $\underline{P}_{j}=P_{0_{j}}\underline{V}_{cl_{j}}^{\alpha_{p_{j}}}$ and $\overline{P}_{j}=P_{0_{j}}\overline{V}_{cl_{j}}^{\alpha_{p_{j}}}$, respectively. To illustrate the relationship between the active power and reactive power of a bus installed with an ES aggregator, an example is shown in Fig. \ref{PQrelationship} with parameters given in Table \ref{table_1}. The blue curves $\underline{g}_{j}(P_{j})$ and $\overline{g}_{j}(P_{j})$ in Fig. \ref{PQrelationship} are the original bounds of the ES aggregator reactive power. For simplicity, the rectangle area in Fig. \ref{PQrelationship} is used as the active and reactive power limits of ES aggregator $j$ where $\underline{Q}_{j}=max~\underline{g}_{j}(P_{j})$ and $\overline{Q}_{j}=min~\overline{g}_{j}(P_{j})$, $\forall j\in\mathcal{N}_{i}\setminus\{1\}$.

\brem
The phase angle differences between each bus, i.e., $|\theta_{j}-\theta_{k}|$, $\forall j,k \in \mathcal{N}_{i}$, in subtransmission networks need to be considered in this optimization problem. Therefore, in each subtransmission network, we assume the phase angle of the PCC as the reference, namely $\theta_1=\theta_{P}$ ($\theta_{P}$ is a constant), in the optimization for simplicity.
\erem

\begin{table}[h]
\renewcommand{\arraystretch}{1}
\centering
\caption{Specifications of the example for ES aggregator operation region}
\label{table_1}
\scalebox{0.8}{
\begin{tabular}{|c|c|}
\hline
Symbol & Value\\
\hline
$V_{s}$ nominal value & 1 p.u.\\
\hline
$V_{es}$ limits & $\pm$1 p.u.\\
\hline
$V_{cl}$ limits & $[0.6,1.4]$ p.u.\\
\hline
$P_{cl}$ nominal value & 1 p.u.\\
\hline
$Q_{cl}$ nominal value & 0.2 p.u.\\
\hline
$\alpha_{p}$ & 1.7\\
\hline
$\alpha_{q}$ & 1.4\\
\hline
\end{tabular}}
\end{table}

To solve this optimization problem in a distributed way, a consensus-based distributed optimization algorithm proposed in \cite{distributedopimalcontinous} is adopted in this paper. It should be noted that matrices $\textbf{B}$, $\textbf{G}$ and $\textbf{G}^{'}$ in (\ref{DLPF2}) are sparse matrices and represent the network topology of the grid. In other words, the injected active power and reactive power of each bus are only related with the voltages of buses to which it has physical connections. Therefore, to obtain the optimal active and reactive power setpoints, the ES aggregator at each bus only needs to calculate its own power flows as well as those of its neighbors it has physical connections to by sharing local information with neighbors. It should be noted that the information of the active power limits of bus $1$, i.e., $\overline{P}_{1}=\underline{P}_{1}=\bar{u}+P_{loss}$, will be broadcast to all ES aggregators periodically, according to which the local constraints of ES aggregators will be updated.

In each subtransmission network, the communication network is considered as an undirected graph, the topology of which is denoted by $C=(c_{jq})\in\mathcal{R}^{N_{i}\times N_{i}}$ where $c_{jj}=0$ and $c_{jq}=c_{qj}>0$ if there is a communication connection between bus $j$ and bus $q$ and $c_{jq}=c_{qj}=0$ otherwise. In this paper, the topology of the communication network is assumed to be the same as that of the physical network since each ES aggregator only needs information from neighbors to implement corresponding power flow calculations. The set of neighbors of ES aggregator $j$ is denoted by $\mathcal{N}_{nb_{j}}=\{w| c_{jw}>0, w\in\mathcal{N}_{i}\}$, and $\mathcal{N}_{j}=\mathcal{N}_{nb_{j}}\cup\{j\}$ and $N_j$ is used to represent the number of elements in each set $\mathcal{N}_{j}$.

Thus, problem (\ref{centraloptimal}) can be transformed to a distributed optimization problem described as the following form,
\begin{equation}\label{distributedoptimal}
\text{minimize}~~f^{total}=\sum^{N^{i}_{es}}_{j=1}f_{j}(\textbf{x})
\end{equation}
subject to
\begin{equation}
\label{distributedoptimal_subject}
\textbf{W}_{j}~\textbf{x}=\boldsymbol{0}_{j}~(j=1,2,...,N_{es}^{i}),~\textbf{x}\in\bigcap^{N_{es}^{i}}_{j=1}\Omega_{j}
\end{equation}
where matrix $\textbf{W}_{j}$ is the submatrix of matrix $\textbf{W}$ defined as
\begin{equation}
\textbf{W}_{j}=\begin{bmatrix}
\textbf{I}_{j} & \textbf{O}_{j} & -\textbf{G}_{j} & \textbf{B}^{'}_{j} \\
\textbf{O}_{j} & \textbf{I}_{j} & \textbf{B}_{j} & \textbf{G}_{j}
\end{bmatrix}
\end{equation}
and the matrices $\textbf{I}_{j},\textbf{O}_{j},\textbf{G}_{j},\textbf{B}_{j}$ and $\textbf{B}^{'}_{j}$ are submatrices which contain the $w$th rows of matrices $\textbf{I},\textbf{O},\textbf{G},\textbf{B}$ and $\textbf{B}^{'}$, respectively, if aggregator $w$ is a neighbor of aggregator $j$, i.e., $w\in\mathcal{N}_{j}$.

Let $\nabla f(x)$ denote the gradient of a function $f(x)$, and $p(v)$ denote a projection operator from $\mathcal{R}^{n}$ to $\Omega\subseteq\mathcal{R}^{n}:~p(v)=argmin_{u\in\Omega}||u-v||$. Let $\textbf{x}_{j}\in\mathcal{R}^{4N_{i}}$ denote an estimated solution to problem (\ref{distributedoptimal}) by ES aggregator $j$. To solve problem (\ref{distributedoptimal}) coorperatively, each ES aggregator $j$ will optimize its local objective function, i.e. $f_{j}(\textbf{x}_{j})$, subject to local constraints, i.e., $\textbf{W}_{j}\textbf{x}_{j}=\boldsymbol{0}_{j}$ and $\textbf{x}_{j}\in\Omega_{j}$, and meanwhile share information $\textbf{x}_{j}$ with neighbors. Thus, each ES aggregator $j$ will generate $\textbf{w}_{j},\textbf{x}_{j},\textbf{y}_{j}$ and $\textbf{z}_{j}$ according to the following rules,

\begin{subequations}\label{iteration}
\begin{align}
\dot{\mathbf{w}}_{j}&=\kappa(-\textbf{w}_{j}\!+\!\textbf{x}_{j}\!-\!\nabla f_{j}(\textbf{x}_{j})\!-\!\textbf{W}_{j}^{T}\textbf{y}_{j}\!-\!\sum^{N_{i}}_{q=1}c_{jq}(\textbf{x}_{j}\!-\!\textbf{x}_{q})\!-\!\textbf{z}_{j})\\
\dot{\mathbf{y}}_{j}&=\zeta\textbf{W}_{j}\textbf{x}_{j}\\
\dot{\mathbf{z}}_{j}&=\eta\sum^{N_{i}}_{q=1}c_{jq}(\textbf{x}_{j}-\textbf{x}_{q})\\
\textbf{x}_{j}&=p^{j}(\textbf{w}_{j})
\end{align}
\end{subequations}
where the vectors $\textbf{w}_{j}\in\mathcal{R}^{4N_{i}}$, $\textbf{y}_{j}\in\mathcal{R}^{2N_{j}}$ and $\textbf{z}_{j}\in\mathcal{R}^{4N_{i}}$ are ancillary vectors, and symbols $\kappa,\zeta$ and $\eta$ are control gains to be designed. The consensus of $\textbf{x}_{j}$ can be achieved consequently from (\ref{iteration}) and meanwhile the optimal solutions of (\ref{distributedoptimal}) can be obtained. The acquired elements $P_{j}$ and $Q_{j}$ in $\textbf{x}_{j}$ are the optimal active and reactive power setpoints of each ES aggregator $j$. 

\section{Case Study}\label{case study}

\begin{figure}[htbp]
\centering
\includegraphics[width=3.5in]{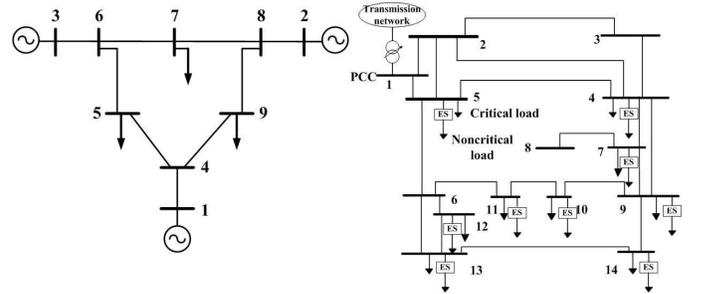}
\caption{IEEE 9-bus system (on the left) and the modified IEEE 14-bus system (on the right)}
\label{case9_14_topo}
\end{figure}

\begin{figure}[htbp]
\centering
\includegraphics[width=2.3in]{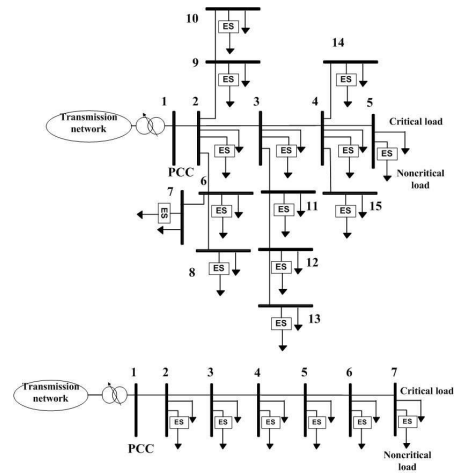}
\caption{The modified IEEE 15-bus system (on the top) and 7-bus feeder system (on the bottom)}
\label{case15_7_topo}
\end{figure}

%
%

%
%
In this section, the IEEE 9-bus system consisting of 3 generators and 3 loads (see Fig. \ref{case9_14_topo}) is used as the transmission network to test the proposed control method. Aggregate loads at bus $5$, $7$ and $9$ are considered as subtransmission networks, respectively. The three subtransmission networks are selected as a 7-bus feeder system, a modified IEEE 15-bus system and a modified IEEE 14-bus system, the configurations of which are shown in Fig. \ref{case15_7_topo} and Fig. \ref{case9_14_topo}, respectively. 
The case with a contingency of a sudden load increase is studied. Comparisons between the proposed control method and the traditional AGC are made in this section.

The nominal value of the frequency and per units of the power and bus voltage are $60$ Hz, $100$ MVA and $220$ KV in the case study, respectively. The nominal value of voltages in the transmission network and subtransmission network are $345$ KV and $220$ KV, respectively. 


At each load bus in the 9-bus system, we assume the load has $20\%$ acceptable adjustment range, i.e., $-0.1P_{D_{i}}\le u_{i}\le 0.1P_{D_{i}}$, and hence the ratio of noncritical load to critical load at each bus in each subtransmission network is $1$ to $9$. For generators in the 9-bus system, we adopt parameters of the governor and reheat steam turbine used in \cite{kundur1994power}. For the selection of parameters of each load-side controller in the 9-bus system, please refer to our previous work \cite{ZHANG_GM} for details. 

It should be noted that all generators in the original 14-bus system and 15-bus system are removed and only loads are preserved in the case study. 
The reactance $X$ and resistance $R$ of all transmission lines in the 7-bus system, 15-bus system and 14-bus system are uniformly set to $0.37$ p.u. and $0.1$ p.u., respectively. 
Nominal active and reactive power consumption of each load bus in each subtransmission network are given in Table \ref{loadpara} in Supporting Document 
which are denoted by $P_n$ and $Q_n$, respectively. For the PCC (bus 1) in each subtransmission network, the bus voltage and reactive power limits are set as $V_{P}=1.05,\theta_{P}=0^{\circ},\underline{Q}_{P}=-100$ and $\overline{Q}_{P}=100$, respectively \cite{zimmerman2011matpower}. For buses $j\in\mathcal{N}_{i}\setminus\{1\}$, we set the limits of magnitudes and phase angles of bus voltages $[\underline{V}_{s},\overline{V}_{s}]=[0.95,1.05]$ p.u. and $[\underline{\theta},\overline{\theta}]=[-15^{\circ},15^{\circ}]$, respectively. Saturation limits of $d$ and $q$ components of the ES voltage setpoints $V_{es\_d_{j}}^{*}$ and $V_{es\_q_{j}}^{*}$ are all set to $-0.7$ p.u. and $0.7$ p.u.. The voltage limits of the noncritical load in each ES aggregator $j$ $[\underline{V}_{cl_{j}},\overline{V}_{cl_{j}}]$ are uniformly set to $[0.6,1.4]$ p.u.. The local constraint of each ES aggregator $j$ denoted by $\Omega_{j}=[\underline{X}_{j},\overline{X}_{j}]$ is uniformly set to be the same as $\Omega$, i.e. $\Omega_{1}=\Omega_{2}=...=\Omega_{N_{es}^{i}}=\Omega$.
For simplicity, cost coefficients $h_{j}$ and $g_{j}$ of each ES aggregator $j$ are uniformly set to $100$ and $40$, respectively. Parameters $\alpha_{p_{j}}$ and $\alpha_{q_{j}}$ of the noncritical load in each ES aggregator $j$ are given in Table \ref{loadpara} in Supporting Document.
 Control gains $\zeta$ and $\eta$ for ES aggregators in 7-bus system, 15-bus system and 14-bus system are uniformly set to $500$ and $250$, respectively, and control gain $\kappa$ is set to $250,190$ and $250$ in each system, respectively. The required active power response is broadcast to aggregators every $0.15$ s within which the new active and reactive power setpoints can be obtained by each aggregator. The linear regression approach is used to acquire the parameters $D^{l}$ and $C$ of each subtransmission network. As a result, $D^{l}$ of the 7-bus system, 15-bus system and 14-bus system are set to $0.0704$, $0.0688$ and $0.0751$, and $C$ are set to $0.0402$, $0.0697$ and $0.0645$, respectively.

\begin{figure}[htbp]
\vspace{-3pt}
\centering
\includegraphics[width=3.0in,height=1.8in]{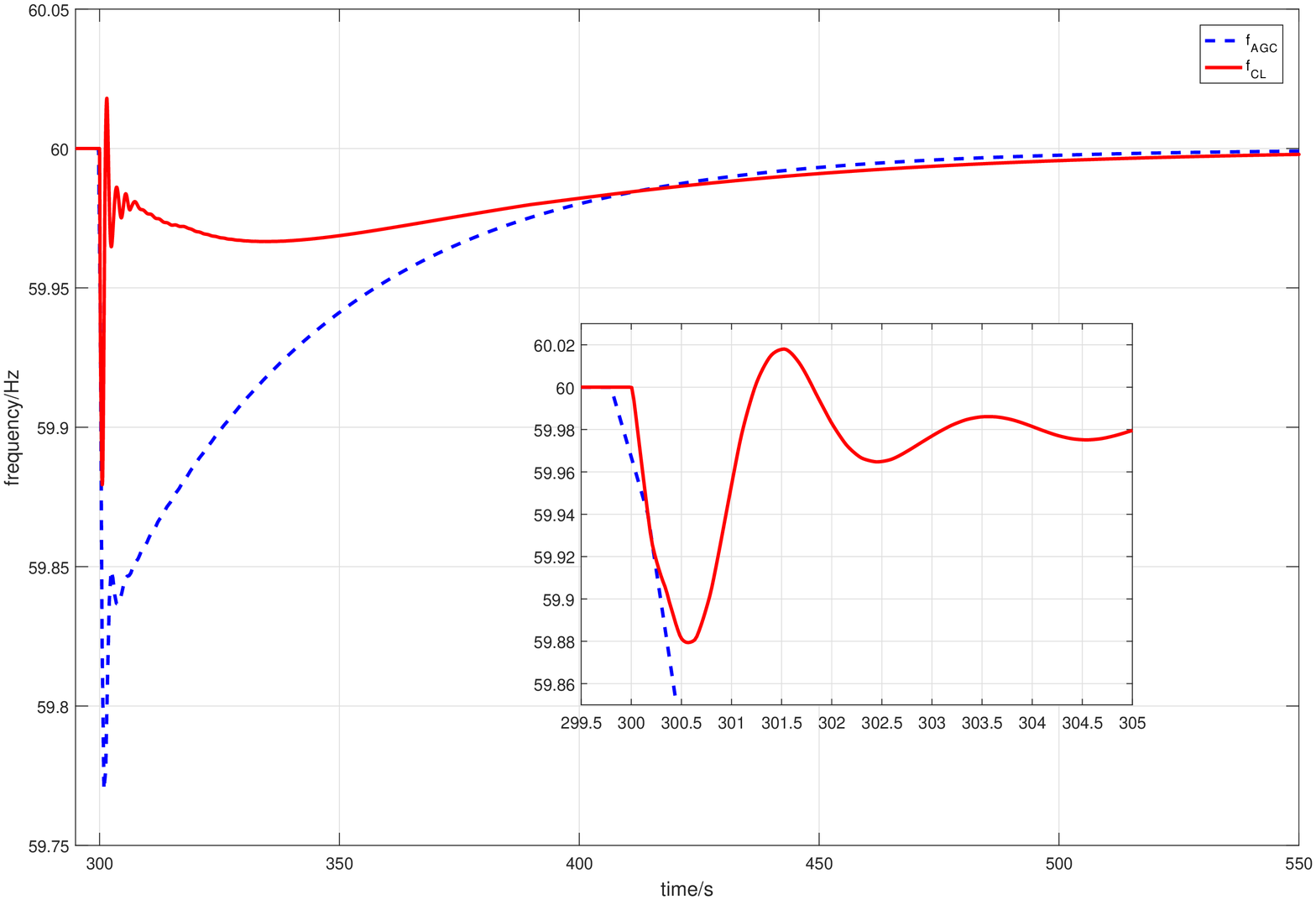}
\caption{The system frequency}
\label{case9_f}
\vspace{-6pt}
\end{figure}

\begin{figure}[htbp]
\vspace{-6pt}
\centering
\includegraphics[width=3.0in,height=1.5in]{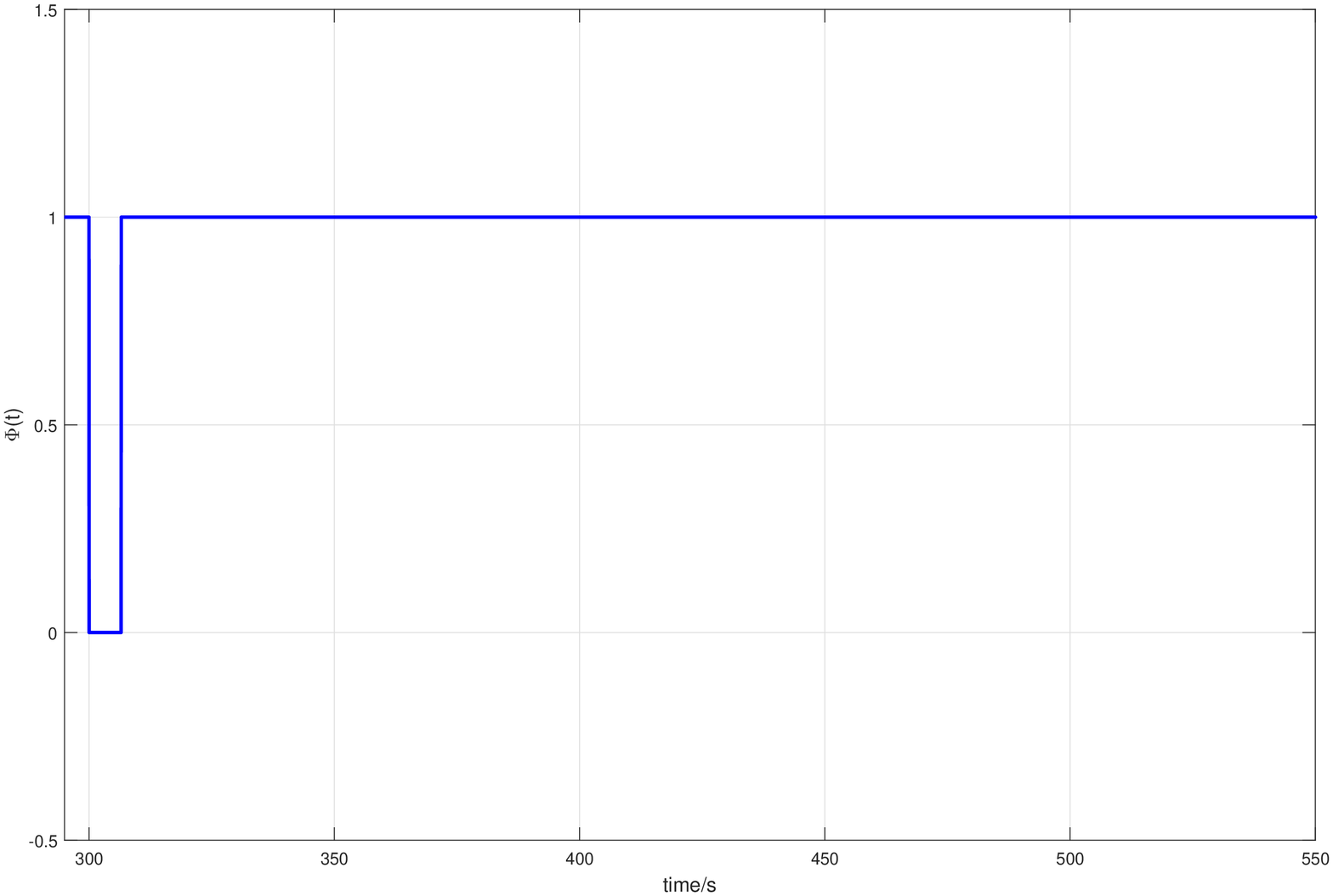}
\caption{The switching signal}
\label{case9_switch}
\vspace{-6pt}
\end{figure}

\begin{figure}[htbp]
\vspace{-6pt}
\centering
\includegraphics[width=3.0in,height=2in]{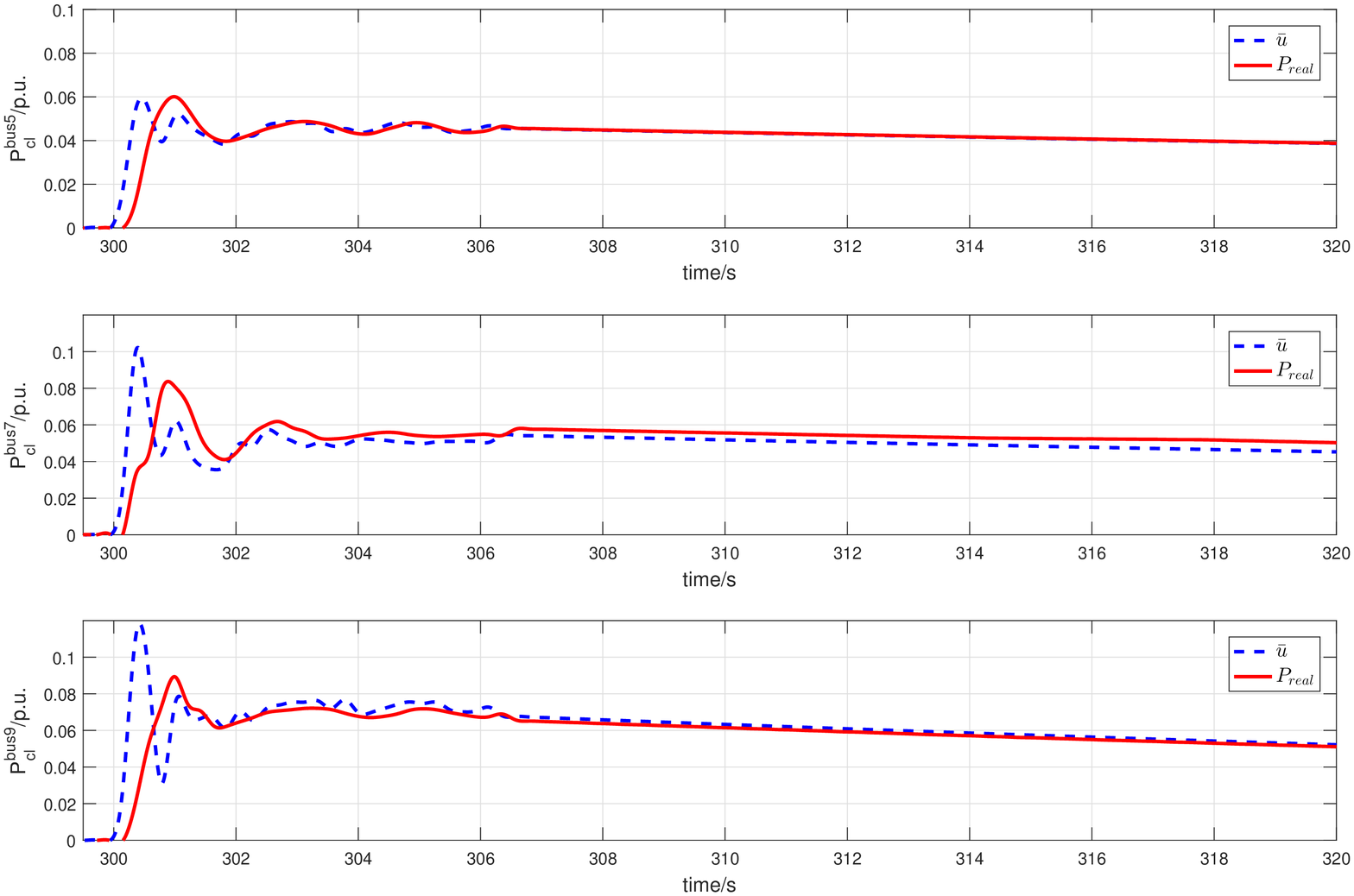}
\caption{Responses of controllable loads in the transmission network and actual aggregate responses of ES aggregators of each subtransmission network}
\label{case9_prs}
\vspace{-6pt}
\end{figure}

\begin{figure}[htbp]
\vspace{-6pt}
\centering
\includegraphics[width=3.0in,height=2in]{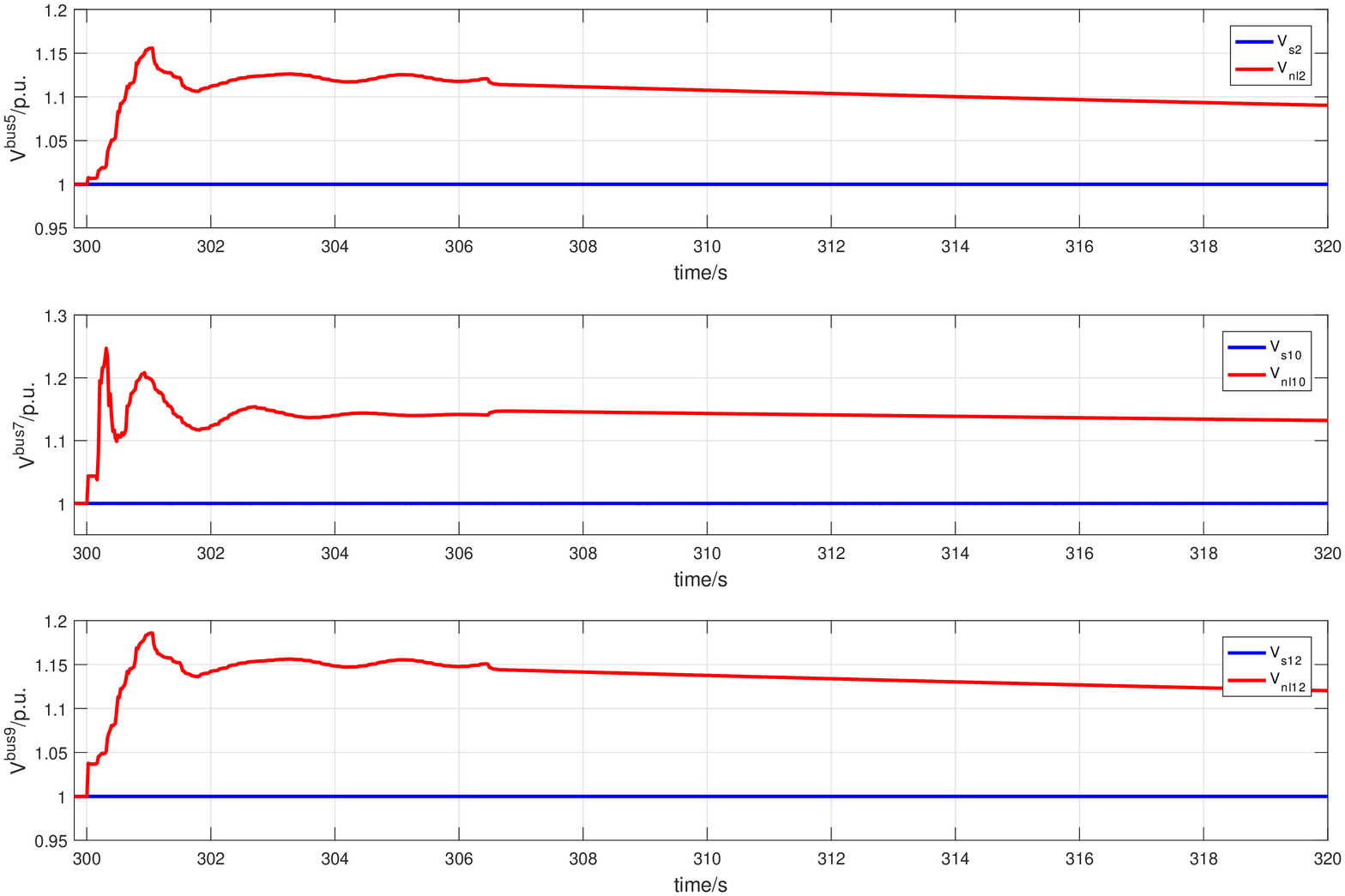}
\caption{Voltages of the critical load and noncritical load at the selected bus in each subtransmission network}
\label{case9_VsVnl}
\vspace{-6pt}
\end{figure}

The system is assumed to operate at its steady state before $t=300$ s, and we assume that a $0.2$ p.u. load increase at bus $7$ in the 15-bus system occurs at $t=300$ s. After the contingency occurs, the control centre detects the system frequency exceeding the frequency band, then it sets $\Phi(t)=0$ to activate the FRM as shown in Fig. \ref{case9_switch}, and then all load-side controllers switch to the FRM immediately. In the meantime, the value of $\bar{u}+P_{loss}$ is sent to each ES aggregator by the control centre, according to which the local constraints $\underline{P}_{1}=\overline{P}_{1}$ of each ES aggregator $j$ are reset. In the meantime, each ES aggregator will share information $\textbf{x}_{j}$ with neighbors to obtain the updated active and reactive power setpoints $P_{j}^{*}$ and $Q_{j}^{*}$ cooperatively, and then adjusts the ES-B2B voltage according to the voltage setpoints obtained in (\ref{apcm_2}). 
 It takes a short time for the distributed optimization to converge to the optimal solution. This leads to a small deviation between the actual aggregate response of ES aggregators in each subtransmission network $P_{real}$ and the required power response $\bar{u}$ as shown in Fig. \ref{case9_prs}. However, it can be observed from Fig. \ref{case9_prs} that $P_{real}$ is still able to track $\bar{u}$ closely during most of the time. As a consequence, the system frequency is restored much more quickly than that under the traditional AGC as shown in Fig. \ref{case9_f}. The voltages of the critical load and noncritical load at bus $2$ in the 7-bus system, bus $12$ in the 15-bus system and bus $10$ in the 14-bus system are shown in Fig. \ref{case9_VsVnl}, respectively. It can be observed from Fig. \ref{case9_VsVnl} that the noncritical load voltages vary greatly when the disturbance occurs to achieve the required active power response and voltage regulation, and consequently the critical voltages only have minor changes. 
As shown in Fig. \ref{case9_switch}, the control centre sets $\Phi(t)=1$ after $t=308$ s when the frequency recovers into the satisfactory frequency region, and then all load-side controllers switch to the LRM. Consequently, power consumption of ES aggregators in each subtransmission network recover to nominal values gradually as shown in Fig. \ref{case9_prs}. Moreover, the costs for ES aggregators power support under the proposed distributed optimization and the proportional adjustment approach (i.e., each ES aggregator adjust its power consumption in proportion to the capacity) are given in Table \ref{cost}, respectively. The percent amounts shown beside the actual difference in Table \ref{cost} correspond to the relative difference over the proportional adjustment approach, which prove the proposed approach is able to reduce costs significantly compared with the ones under the proportional adjustment approach.

\begin{table}[h]
\renewcommand{\arraystretch}{1}
\caption{Costs for ES aggregator active power compensation}
\label{cost}
\scalebox{0.8}{
\begin{tabular}{|c|c|c|c|}
\hline
Network & Costs with PA \footnotemark[1] & Costs with DO \footnotemark[2] & Difference\\
\hline
7-bus system & 11.62 & 3.08 & 8.54 (73.5\%)\\
\hline
15-bus system & 6.33 & 2.81 & 3.52 (55.6\%)\\
\hline
14-bus system & 10.09 & 4.93 & 5.16 (51.1\%)\\
\hline
\end{tabular}}
\footnotemark[1]{PA=proportional adjustment}
\footnotemark[2]{DO=distributed optimization}
\end{table}

\section{Conclusion}

To achieve required performances of load-side controllers in transmission networks and regulate frequency, aggregate controllable loads need to be granulated down to subtransmission networks where both frequency and bus voltage are affected by active power changes due to a higher $R/X$ ratio of transmission lines. Further, the costs for demand response also need to be considered when controllable loads participate in frequency regulation. In this paper a control scheme has been proposed for subtransmission networks in which a distributed optimization is adopted for each electric spring aggregator to obtain the updated active and reactive power setpoints and adjust voltage accordingly. The simulation results have shown that the required active power compensation can be implemented by ES aggregators cooperatively and bus voltages are maintained within the acceptable range all the time with the minimized costs under the proposed control scheme.

\bibliographystyle{IEEEtran}
\bibliography{ZLH_ESoptimal_arXiv}

\begin{thebibliography}{10}
\providecommand{\url}[1]{#1}
\csname url@samestyle\endcsname
\providecommand{\newblock}{\relax}
\providecommand{\bibinfo}[2]{#2}
\providecommand{\BIBentrySTDinterwordspacing}{\spaceskip=0pt\relax}
\providecommand{\BIBentryALTinterwordstretchfactor}{4}
\providecommand{\BIBentryALTinterwordspacing}{\spaceskip=\fontdimen2\font plus
\BIBentryALTinterwordstretchfactor\fontdimen3\font minus
  \fontdimen4\font\relax}
\providecommand{\BIBforeignlanguage}[2]{{%
\expandafter\ifx\csname l@#1\endcsname\relax
\typeout{** WARNING: IEEEtran.bst: No hyphenation pattern has been}%
\typeout{** loaded for the language `#1'. Using the pattern for}%
\typeout{** the default language instead.}%
\else
\language=\csname l@#1\endcsname
\fi
#2}}
\providecommand{\BIBdecl}{\relax}
\BIBdecl

\bibitem{kundur1994power}
P.~Kundur, \emph{Power System Stability and Control}.\hskip 1em plus 0.5em
  minus 0.4em\relax McGraw-hill New York, 1994.

\bibitem{windenergy_2005}
R.~Piwko, D.~Osborn, R.~Gramlich, G.~Jordan, D.~Hawkins, and K.~Porter, ``Wind
  energy delivery issues: transmission planning and competitive electricity
  market operation,'' \emph{Power and Energy Magazine, IEEE}, vol.~3, no.~6,
  pp. 47--56, Nov 2005.

\bibitem{xu2011demand}
Z.~Xu, J.~{\O}stergaard, and M.~Togeby, ``Demand as frequency controlled
  reserve,'' \emph{Power Systems, IEEE Transactions on}, vol.~26, no.~3, pp.
  1062--1071, 2011.

\bibitem{callaway2011achieving}
D.~S. Callaway and I.~Hiskens, ``Achieving controllability of electric loads,''
  \emph{Proceedings of the IEEE}, vol.~99, no.~1, pp. 184--199, 2011.

\bibitem{tao_nondisrupt}
T.~Liu, D.~J. Hill, and C.~Zhang, ``Non-disruptive load-side control for
  frequency regulation in power systems,'' \emph{Smart Grid, IEEE Transactions
  on}, vol.~7, no.~4, pp. 2142--2153, July 2016.

\bibitem{ZHANG_GM}
C.~Zhang, T.~Liu, and D.~J. Hill, ``Switched distributed load-side frequency
  regulation for power systems,'' in \emph{Power and Energy Society General
  Meeting, IEEE}, 2016, pp. 1--5.

\bibitem{zhao2015distributed}
C.~Zhao, E.~Mallada, and F.~Dorfler, ``Distributed frequency control for
  stability and economic dispatch in power networks,'' in \emph{American
  Control Conference, IEEE}, 2015, pp. 2359--2364.

\bibitem{zhangpscc}
C.~Zhang, T.~Liu, and D.~J. Hill, ``Distributed load-side frequency regulation
  for power systems,'' in \emph{Power Systems Computation Conference}, June
  2016, pp. 1--7.

\bibitem{mallada2014distributed}
E.~Mallada and S.~H. Low, ``Distributed frequency-preserving optimal load
  control,'' in \emph{IFAC World Congress}, 2014, pp. 5411--5418.

\bibitem{Zhengyu_GM}
Y.~Zheng, D.~J. Hill, C.~Zhang, and K.~Meng, ``Non-interruptive
  thermostatically controlled load for primary frequency support,'' in
  \emph{Power and Energy Society General Meeting, IEEE}, 2016, pp. 1--5.

\bibitem{hierarchy}
M.~D. Ilic, ``From hierarchical to open access electric power systems,''
  \emph{Proceedings of the IEEE}, vol.~95, no.~5, pp. 1060--1084, May 2007.

\bibitem{zhangccpowertech}
C.~Zhang, T.~Liu, and D.~J. Hill, ``Granulated load-side control of power
  systems with electric spring aggregators,'' in \emph{12th IEEE PES PowerTech
  Conference}, 2017.

\bibitem{b2b_voltageES}
Z.~Akhtar, B.~Chaudhuri, and S.~Y.~R. Hui, ``Smart loads for voltage control in
  distribution networks,'' \emph{IEEE Transactions on Smart Grid}, vol.~PP,
  no.~99, pp. 1--10, 2016.

\bibitem{yan_esb2b}
S.~Yan, C.~K. Lee, T.~B. Yang, K.~T. Mok, S.~C. Tan, B.~Chaudhuri, and S.~Y.~R.
  Hui, ``Extending the operating range of electric spring using back-to-back
  converters: Hardware implementation and control,'' \emph{IEEE Transactions on
  Power Electronics}, vol.~PP, no.~99, pp. 1--1, 2016.

\bibitem{hui2012electric}
S.~Y. Hui, C.~K. Lee, and F.~F. Wu, ``Electric springs---a new smart grid
  technology,'' \emph{Smart Grid, IEEE Transactions on}, vol.~3, no.~3, pp.
  1552--1561, 2012.

\bibitem{DLPF}
J.~Yang, N.~Zhang, C.~Kang, and Q.~Xia, ``A state-independent linear power flow
  model with accurate estimation of voltage magnitude,'' \emph{IEEE
  Transactions on Power Systems}, vol.~PP, no.~99, pp. 1--1, 2016.

\bibitem{bergen1981structure}
A.~R. Bergen and D.~J. Hill, ``A structure preserving model for power system
  stability analysis,'' \emph{Power Apparatus and Systems, IEEE Transactions
  on}, vol.~1, pp. 25--35, 1981.

\bibitem{zheng2017critical}
Y.~Zheng, D.~J. Hill, K.~Meng, and S.~Hui, ``Critical bus voltage support in
  distribution systems with electric springs and responsibility sharing,''
  \emph{IEEE Transactions on Power Systems}, vol.~32, no.~5, pp. 3584--3593,
  2017.

\bibitem{wu2016line}
A.~Wu and B.~Ni, \emph{Line Loss Analysis and Calculation of Electric Power
  Systems}.\hskip 1em plus 0.5em minus 0.4em\relax John Wiley \& Sons, 2016.

\bibitem{distributedopimalcontinous}
S.~Yang, Q.~Liu, and J.~Wang, ``A multi-agent system with a
  proportional-integral protocol for distributed constrained optimization,''
  \emph{IEEE Transactions on Automatic Control}, vol.~62, no.~7, pp.
  3461--3467, July 2017.

\bibitem{zimmerman2011matpower}
R.~D. Zimmerman, C.~E. Murillo-S{\'a}nchez, and R.~J. Thomas, ``Matpower:
  Steady-state operations, planning, and analysis tools for power systems
  research and education,'' \emph{Power Systems, IEEE Transactions on},
  vol.~26, no.~1, pp. 12--19, 2011.

\end{thebibliography}

\clearpage

\section{Supporting Document}

\emph{Derivation of (\ref{apcm_2}):} Subtracting (\ref{Qnl}) from (\ref{Pnl}) leads to
 \begin{equation}\label{Vnldq}
 (V_{nl\_d_{j}}^{2}+V_{nl\_q_{j}}^{2})I_{q_{j}}=P_{nl_{j}}V_{nl\_q_{j}}-Q_{nl_{j}}V_{nl\_d_{j}}.
\end{equation}
By substituting (\ref{Vnl}), (\ref{sqrt}) and (\ref{Pnl}) into (\ref{Vnldq}), we can acquire
\begin{equation}\label{VsVesd}
V_{nl\_d_{j}}=-\frac{V_{0_{j}}^{2}\!\left(\!\frac{P^{*}_{j}}{P_{0_{j}}}\!\right)\!^{\frac{2}{\alpha_{p_{j}}}}I_{q_{j}}+P^{*}_{j}V_{es\_q_{j}}}{Q_{nl_{j}}}.
\end{equation}
Substituting (\ref{VsVesd}) into (\ref{sqrt}) and replacing $V_{nl\_q_{j}}$ by $-V_{es\_q_{j}}$ lead to
\begin{equation}\label{Vesq2}
\begin{split}
(P^{*2}_{j}+Q^{2}_{nl_{j}})V_{es\_q_{j}}^{2}+2V^{2}_{0_{j}}\!\left(\!\frac{P^{*}_{j}}{P_{0_{j}}}\!\right)\!^{\frac{2}{\alpha_{p_{j}}}}I_{q_{j}}P^{*}_{j}V_{es\_q_{j}}\\
+V^{4}_{0_{j}}\!\left(\!\frac{P^{*}_{j}}{P_{0_{j}}}\!\right)\!^{\frac{4}{\alpha_{p_{j}}}}I_{q_{j}}^{2}-V^{2}_{0_{j}}\!\left(\!\frac{P^{*}_{j}}{P_{0_{j}}}\!\right)\!^{\frac{2}{\alpha_{p_{j}}}}Q^{2}_{nl_{j}}=0.
\end{split}
\end{equation}
By substituting (\ref{Qnl}) into (\ref{Vesq2}) and using $I_{q_{j}}=-\frac{Q^{*}_{j}}{V_{s_{j}}}$ we have
\begin{equation}\label{Vesq3}
\begin{split}
\!\left(\!P^{*2}_{j}+Q_{0_{j}}^{2}\!\left(\!\frac{P^{*}_{j}}{P_{0_{j}}}\!\right)\!^{\frac{2\alpha_{q_{j}}}{\alpha_{p_{j}}}}\!\right)\!V_{es\_q_{j}}^{2}\!-\!2V^{2}_{0_{j}}\!\left(\!\frac{P^{*}_{j}}{P_{0_{j}}}\!\right)\!^{\frac{2}{\alpha_{p_{j}}}}\frac{Q^{*}_{j}}{V_{s_{j}}}P^{*}_{j}V_{es\_q_{j}}\\
\!+\!\!\left(\!V^{4}_{0_{j}}\!\left(\!\frac{P^{*}_{j}}{P_{0_{j}}}\!\right)\!^{\frac{4}{\alpha_{p_{j}}}}\frac{Q^{*2}_{j}}{V^{2}_{s_{j}}}\!-\!V^{2}_{0_{j}}\!\left(\!\frac{P^{*}_{j}}{P_{0_{j}}}\!\right)\!^{\frac{2}{\alpha_{p_{j}}}}Q_{0_{j}}^{2}\!\left(\!\frac{P^{*}_{j}}{P_{0_{j}}}\!\right)\!^{\frac{2\alpha_{q_{j}}}{\alpha_{p_{j}}}}\!\right)\!\!=\!0
\end{split}
\end{equation}
For simplicity, here we assume that the bus voltage can be maintained at the setpoint, i.e., $V_{s_{j}}=V^{*}_{s_{j}}$, due to the control actions of ES aggregators. It is reasonable to set $V^{*}_{s_{j}}=V_{0_{j}}$, and hence we have $V_{s_{j}}=V^{*}_{s_{j}}=V_{0_{j}}$ and (\ref{Vesq3}) can be rewritten as follows,

\begin{equation}\label{Vesq4}
\begin{split}
\!\left(\!P^{*2}_{j}+Q_{0_{j}}^{2}\!\left(\!\frac{P^{*}_{j}}{P_{0_{j}}}\!\right)\!^{\frac{2\alpha_{q_{j}}}{\alpha_{p_{j}}}}\!\right)\!V_{es\_q_{j}}^{2}\!-\!2V_{0_{j}}\!\left(\!\frac{P^{*}_{j}}{P_{0_{j}}}\!\right)\!^{\frac{2}{\alpha_{p_{j}}}}Q^{*}_{j}P^{*}_{j}V_{es\_q_{j}}\\
\!+\!\!\left(\!V^{2}_{0_{j}}\!\left(\!\frac{P^{*}_{j}}{P_{0_{j}}}\!\right)\!^{\frac{4}{\alpha_{p_{j}}}}Q^{*2}_{j}\!-\!V^{2}_{0_{j}}\!\left(\!\frac{P^{*}_{j}}{P_{0_{j}}}\!\right)\!^{\frac{2}{\alpha_{p_{j}}}}Q_{0_{j}}^{2}\!\left(\!\frac{P^{*}_{j}}{P_{0_{j}}}\!\right)\!^{\frac{2\alpha_{q_{j}}}{\alpha_{p_{j}}}}\!\right)\!\!=\!0
\end{split}
\end{equation}

By solving (\ref{Vesq4}), we can acquire the $q$ component of the ES-B2B voltage setpoint $V^{*}_{es\_q_{j}}$ as follows,
\begin{equation}\label{Vesq3}
V^{*}_{es\_q_{j}}=\frac{\left(\frac{P^{*}_{j}}{P_{0_{j}}}\right)^{\frac{2}{\alpha_{p_{j}}}}P^{*}_{j}Q^{*}_{j}-Q_{0_{j}}\left(\frac{P^{*}_{j}}{P_{0_{j}}}\right)^{\frac{\alpha_{q_{j}}+1}{\alpha_{p_{j}}}}\sqrt{\Delta}}{P^{*2}_{j}+Q_{0_{j}}^{2}\left(\frac{P^{*}_{j}}{P_{0_{j}}}\right)^{\frac{2\alpha_{q_{j}}}{\alpha_{p_{j}}}}}V_{0_{j}},
\end{equation}
where $\Delta=P^{*2}_{j}-\left(\frac{P^{*}_{j}}{P_{0_{j}}}\right)^{\frac{2}{\alpha_{p_{j}}}}Q^{*2}_{j}+Q_{0_{j}}^{2}\left(\frac{P^{*}_{j}}{P_{0_{j}}}\right)^{\frac{2\alpha_{q_{j}}}{\alpha_{p_{j}}}}$. The $d$ component of the ES-B2B voltage setpoint $V^{*}_{es\_d_{j}}$ can be obtained by substituting (\ref{Vesq3}) into (\ref{VsVesd}) as follows,
\begin{equation}
V^{*}_{es\_d_{j}}\!=\!V^{*}_{s_{j}}\!-\!\frac{Q_{0_{j}}\left(\frac{P^{*}_{j}}{P_{0_{j}}}\right)^{\frac{\alpha_{q_{j}}+2}{\alpha_{p_{j}}}}Q^{*}_{j}+\left(\frac{P^{*}_{j}}{P_{0_{j}}}\right)^{\frac{1}{\alpha_{p_{j}}}}P^{*}_{j}\sqrt{\Delta}}{P^{*2}_{j}+Q_{0_{j}}^{2}\left(\frac{P^{*}_{j}}{P_{0_{j}}}\right)^{\frac{2\alpha_{q_{j}}}{\alpha_{p_{j}}}}}V_{0_{j}}.
\end{equation}

\begin{table}[h]
\renewcommand{\arraystretch}{1}
\caption{Load parameters of the noncritical load in each ES aggregator}
\label{loadpara}
\scalebox{0.8}{
\begin{tabular}{|c|c|c|c|c|c|} 
\hline
\multirow{2}{*}{Subtransmission network} & \multirow{2}{*}{Bus No.} & \multicolumn{2}{c|}{Load parameters}& \multirow{2}{*}{$P_n$ (p.u.)} & \multirow{2}{*}{$Q_n$ (p.u.)}\\
\cline{3-4}
 & & $\alpha_{p}$ & $\alpha_{q}$ & & \\
 \hline
 \multirow{6}{*}{7-bus system}
 &2& 1.3 & 1.2 & 0.25 & 0.05\\
 \cline{2-6}
 &3& 1.4 & 1.3 & 0.1 & 0.02\\
  \cline{2-6}
 &4& 1.5 & 1.4 & 0.025 & 0.005\\
  \cline{2-6}
 &5& 1.6 & 1.5 & 0.025 & 0.005\\
  \cline{2-6}
 &6& 1.7 & 1.6 & 0.025 & 0.005\\
  \cline{2-6}
 &7& 1.8 & 1.7 & 0.025 & 0.005\\
  \cline{2-6}
 \hline
  \multirow{14}{*}{15-bus system}
 &2& 1.9 & 1.9 & 0.017 & 0.0035\\
 \cline{2-6}

 &3& 1.9 & 1.9 & 0.0269 & 0.0055\\
  \cline{2-6}

 &4& 1.9 & 1.9 & 0.0538 & 0.0011\\
  \cline{2-6}
 
 &5& 1.8 & 1.8 & 0.017 & 0.0035\\
  \cline{2-6}
 
 &6& 1.7 & 1.7 & 0.0538 & 0.0011\\
  \cline{2-6}
 
 &7& 1.7 & 1.7 & 0.0538 & 0.0011\\
  \cline{2-6}
 
 &8& 1.6 & 1.6 & 0.0269 & 0.0055\\
 \cline{2-6}
 
 &9& 1.6 & 1.6 & 0.0269 & 0.0055\\
  \cline{2-6}
 
 &10& 1.5 & 1.5 & 0.017 & 0.0035\\
  \cline{2-6}
 
 &11& 1.4 & 1.4 & 0.0538 & 0.0011\\
  \cline{2-6}
 
 &12& 1.3 & 1.3 & 0.0269 & 0.0055\\
  \cline{2-6}
 
 &13& 1.2 & 1.2 & 0.017 & 0.0035\\
   \cline{2-6}
 
 &14& 1.1 & 1.1 & 0.0538 & 0.0011\\
  \cline{2-6}
 
 &15& 1.1 & 1.1 & 0.0538 & 0.0011\\
 \hline
  \multirow{9}{*}{14-bus system}
 &4& 1.9 & 1.9 & 0.239 & 0.0195\\
 \cline{2-6}
 &5& 1.8 & 1.8 & 0.038 & 0.008\\
  \cline{2-6}
 &7& 1.7 & 1.7 & 0.056 & 0.0001\\
  \cline{2-6}
 &9& 1.6 & 1.6 & 0.1475 & 0.083\\
  \cline{2-6}
 &10& 1.5 & 1.5 & 0.045 & 0.029\\
  \cline{2-6}
 &11& 1.4 & 1.4 & 0.0175 & 0.009\\
  \cline{2-6}
   &12& 1.3 & 1.3 & 0.0305 & 0.008\\
  \cline{2-6}
 &13& 1.2 & 1.2 & 0.0675 & 0.029\\
  \cline{2-6}
 &14& 1.1 & 1.1 & 0.0745 & 0.025\\
 \hline
 \end{tabular}}
 \end{table}

\end{document}